# DIRICHLET PROBLEM IN THE NON-CLASSICAL TREATMENT FOR ONE PSEUDOPARABOLIC EQUATION OF FOURTH ORDER

**I.G. Mamedov**

*A.I. Huseynov Institute of Cybernetics of NAS of Azerbaijan. Az 1141, Azerbaijan, Baku st. B. Vahabzade, 9*
*E-mail: ilgar-mammadov@rambler.ru*

**Abstract**

*In the paper the Dirichlet problem with non-classical conditions not requiring agreement conditions is considered for a fourth order pseudoparabolic equation with non-classical coefficients. The equivalence of these conditions with the classic boundary conditions is substantiated in the case if the solution of the stated problem is sought in S.L. Sobolev isotropic space.*

## Problem Statement

The first boundary value problem or the Dirichlet problem (i.e. a problem in which a closed contour is an input medium) known well for elliptic type differential equations is one of the basic boundary value problems of mathematical physics [1-2]. From this point of view this paper is devoted to urgent problems of mathematical physics.

Consider the equation

$$(V_{2,2}\,u)(x) \equiv D_1^2 D_2^2 u(x) + a_{2,1}(x) D_1^2 D_2 u(x) + a_{1,2}(x) D_1 D_2^2 u(x) + \\ + a_{2,0}(x) D_1^2 u(x) + a_{0,2}(x) D_2^2 u(x) + \sum_{i_1=0}^{1}\sum_{i_2=0}^{1} a_{i_1,i_2}(x) D_1^{i_1} D_2^{i_2} u(x) = Z_{2,2}(x) \in L_p(G). \quad (1)$$

Here $u(x) \equiv u(x_1, x_2)$ is a desired function determined on $G$; $a_{i_1,i_2}(x)$ are the given measurable functions on $G = G_1 \times G_2$, where $G_j = (0, h_j)$, $j=1,2$; $Z_{2,2}(x)$ is a given measurable function on $G$; $D_j \equiv \partial/\partial x_j$ is a generalized differentiation operator in S.L. Sobolev sense, $j = 1,2$.

Equation (1) is a hyperbolic equation that possesses two real characteristics $x_1 = const$, $x_2 = const$, the first and the second one of which are two-fold. The equation of type (1) in the paper of A.P.Soldatov and M. Kh.Shkhanukov [3] are called pseudoparabolic ones. Different special cases of equation (1) arise by modeling various processes of applied character (generalized equation of moisture transfer, telegraph equation, string vibration equation, heat conductivity equation, Aller's equations and etc.). Furthermore, this equation is a generalization of the Boussenesq-Liav equation [4] describing longitudinal waves in a thin elastic bar with regard to lateral inertia effects.

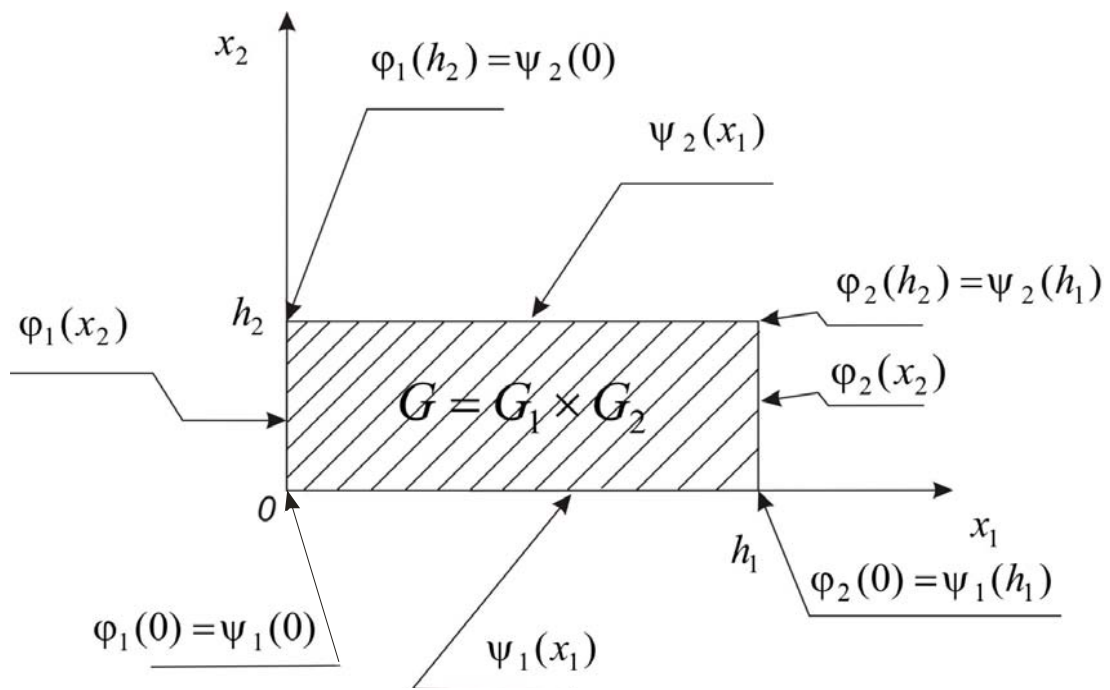

**Fig. 1. Geometrical interpretation of Dirichlet classic conditions**

In the paper we consider equation (1) in the general case when the coefficients $a_{i_1,i_2}(x)$ are non-smooth functions satisfying only the following conditions:

$$a_{2,i_2}(x) \in L^{x_1,x_2}_{\infty,p}(G), \ i_2 = \overline{0,1}; a_{i_1,2}(x) \in L^{x_1,x_2}_{p,\infty}(G), \ i_1 = \overline{0,1};$$

$$a_{i_1,i_2}(x) \in L_p(G), \ i_1 = \overline{0,1}, \ i_2 = \overline{0,1}.$$

Therewith the important principal moment is that the equation under consideration possesses non-smooth coefficients that satisfy only some $p$-integrability and boundedness conditions, i.e. the considered pseudoparabolic differential operator $V_{2,2}$ has no traditional conjugation operator. In other words, the Green function - the source function for such an equation can't be investigated by the classic method of characteristics.

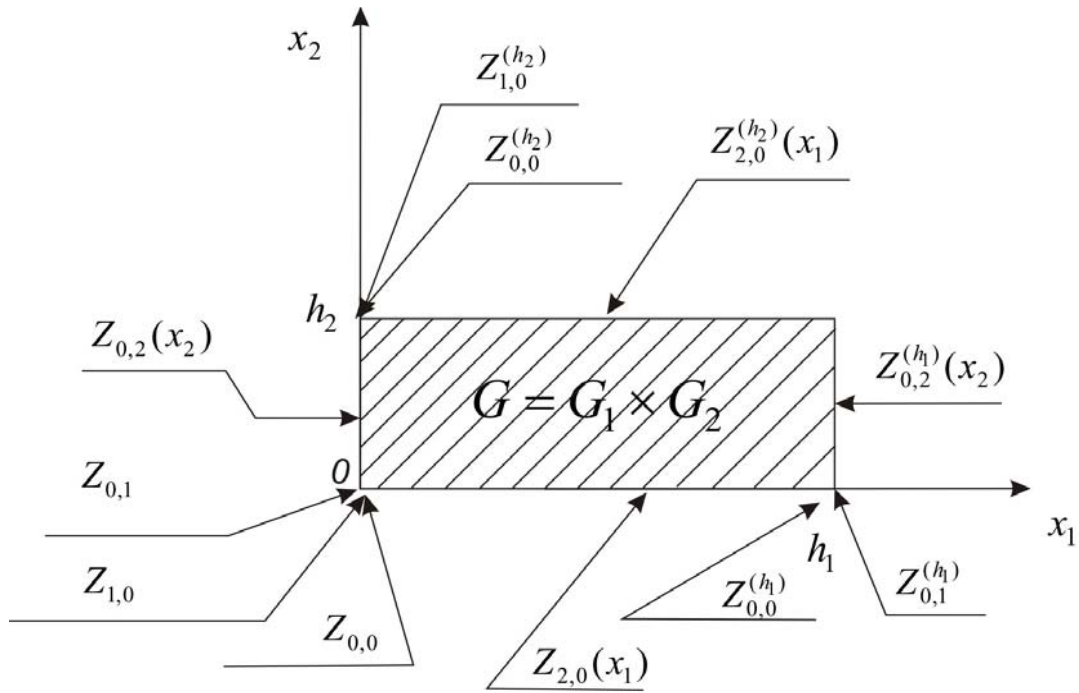

**Fig. 2. Geometrical interpretation of Dirichlet non-classical conditions.**

Under these conditions, we'll look for the solution $u(x)$ of equation (1) in S.L.Sobolev space

$$W_p^{(2,2)}(G) \equiv \left\{ u(x) : D_1^{i_1} D_2^{i_2} u(x) \in L_p(G); \; i_j = \overline{0,2}, \; j = \overline{0,1} \right\},$$

where $1 \leq p \leq \infty$. We'll define the norm in space $W_p^{(2,2)}(G)$ by the equality:

$$\|u(x)\|_{W_p^{(2,2)}(G)} = \sum_{i_1=0}^{2} \sum_{i_2=0}^{2} \left\| D_1^{i_1} D_2^{i_2} u(x) \right\|_{L_p(G)}.$$

For equation (1) we can give the classic form Dirichlet condition [5] as follows

$$\begin{cases} u(0, x_2) = \varphi_1(x_2), \ u(x_1, 0) = \psi_1(x_1), \\ u(h_1, x_2) = \varphi_2(x_2), \ u(x_1, h_2) = \psi_2(x_1), \end{cases} \quad (2)$$

where $\varphi_j(x_2)$ and $\psi_j(x_1)$, $j = \overline{1,2}$ are the given measurable functions on $G$.

Obviously, in the case of conditions (2), in addition to conditions

$$\varphi_j(x_2) \in W_p^{(2)}(G_2) \equiv \{\tilde{\varphi}(x_2): D_2^{i_2}\tilde{\varphi}(x_2) \in L_p(G_2), i_2 = \overline{0,2}\}, 1 \le p \le \infty, j = \overline{1,2};$$
$$\psi_j(x_1) \in W_p^{(2)}(G_1) \equiv \{\tilde{\psi}(x_1): D_1^{i_1}\tilde{\psi}(x_1) \in L_p(G_1), i_1 = \overline{0,2}\}, 1 \le p \le \infty, j = \overline{1,2};$$

the given functions should also satisfy the following agreement conditions:

$$\begin{cases} \varphi_1(0) = \psi_1(0), \ \varphi_2(h_2) = \psi_2(h_1), \\ \varphi_1(h_2) = \psi_2(0), \ \varphi_2(0) = \psi_1(h_1). \end{cases} \quad (3)$$

Consider the following non-classical boundary conditions:

$$\begin{cases} V_{0,0}u \equiv u(0,0) = Z_{0,0} \in R; \\ V_{1,0}u \equiv D_1 u(0,0) = Z_{1,0} \in R; \\ V_{0,1}u \equiv D_2 u(0,0) = Z_{0,1} \in R; \\ (V_{2,0}u)(x_1) \equiv D_1^2 u(x_1, 0) = Z_{2,0}(x_1) \in L_p(G_1); \\ (V_{0,2}u)(x_2) \equiv D_2^2 u(0, x_2) = Z_{0,2}(x_2) \in L_p(G_2); \\ V_{0,0}^{(h_1)}u \equiv u(h_1, 0) = Z_{0,0}^{(h_1)} \in R; \\ V_{0,1}^{(h_1)}u \equiv D_2 u(h_1, 0) = Z_{0,1}^{(h_1)} \in R; \\ V_{0,0}^{(h_2)}u \equiv u(0, h_2) = Z_{0,0}^{(h_2)} \in R; \\ V_{1,0}^{(h_2)}u \equiv D_1 u(0, h_2) = Z_{1,0}^{(h_2)} \in R; \\ (V_{2,0}^{(h_2)}u)(x_1) \equiv D_1^2 u(x_1, h_2) = Z_{2,0}^{(h_2)}(x_1) \in L_p(G_1); \\ (V_{0,2}^{(h_1)}u)(x_2) \equiv D_2^2 u(h_1, x_2) = Z_{0,2}^{(h_1)}(x_2) \in L_p(G_2); \end{cases} \quad (4)$$

If the function $u \in W_p^{(2,2)}(G)$ is a solution of the classic form Dirichlet problem (1), (2), then it is a solution of problem (1), (4) for $Z_{i_1,i_2}$ and $Z_{i_1,i_2}^{(h_j)}$, defined by the following equalities:

$$Z_{0,0} = \varphi_1(0) = \psi_1(0); \ Z_{1,0} = \psi_1'(0); \ Z_{0,1} = \varphi_1'(0);$$
$$Z_{2,0}(x_1) = \psi_1''(x_1); \ Z_{0,2}(x_2) = \varphi_1''(x_2); \ Z_{0,0}^{(h_1)} = \varphi_2(0) = \psi_1(h_1); \ Z_{0,1}^{(h_1)} = \varphi_2'(0);$$
$$Z_{0,0}^{(h_2)} = \psi_2(0) = \varphi_1(h_2); \ Z_{1,0}^{(h_2)} = \psi_2'(0); \ Z_{2,0}^{(h_2)}(x_1) = \psi_2''(x_1); \ Z_{0,2}^{(h_1)}(x_2) = \varphi_2''(x_2).$$

It is easy to prove that inverse one is also true In other words, if the function $u \in W_p^{(2,2)}(G)$ is a solution of problem (1), (4) then it is also a solution of problem (1), (2) for the following functions:

$$\varphi_1(x_2) = Z_{0,0} + x_2 Z_{0,1} + \int_0^{x_2} (x_2 - \tau) Z_{0,2}(\tau) d\tau; \qquad (5)$$

$$\varphi_2(x_2) = Z_{0,0}^{(h_1)} + x_2 Z_{0,1}^{(h_1)} + \int_0^{x_2} (x_2 - \xi) Z_{0,2}^{(h_1)}(\xi) d\xi; \qquad (6)$$

$$\psi_1(x_1) = Z_{0,0} + x_1 Z_{1,0} + \int_0^{x_1} (x_1 - \eta) Z_{2,0}(\eta) d\eta; \qquad (7)$$

$$\psi_2(x_1) = Z_{0,0}^{(h_2)} + x_1 Z_{1,0}^{(h_2)} + \int_0^{x_1} (x_1 - \nu) Z_{2,0}^{(h_2)}(\nu) d\nu; \qquad (8)$$

Note that the functions (5) - (8) possess an important property, more exactly, agreement conditions (3) for all $Z_{i_1,i_2}$ and $Z_{i_1,i_2}^{(h_j)}$, possessing the above-mentioned properties are fulfilled for them automatically. Therefore, we can consider equalities (5)-(8) as a general form of all the functions $\varphi_j(x_2), \psi_j(x_1), j = \overline{1,2}$, satisfying the agreement conditions (3).

So, the classic form Diriclet problems (1), (2) and of the form (1), (4) are equivalent in the general case. However, the Dirichlet problem (1), (4) in nonclassic treatment is more natural by statement than Dirichlet problem (1), (2). This is connected with the fact that in the statement of Dirichlet problem (1), (4), the right sides of boundary conditions have no additional conditions of agreement type.

It should be especially noted that in the papers [6-7] the author suggested a method for investigating boundary value problems in non-classical treatment for pseudoparabolic equations with non-smooth coefficients of higher order.

Note that the Dirichlet non-classical problem in treatment (1), (4) is investigated by means of integral representations of special type for the functions $u(x) \in W_p^{(2,2)}(G)$:

$$u(x) = u(0,0) + x_1 D_1 u(0,0) + x_2 D_2 u(0,0) + x_1 x_2 D_1 D_2 u(0,0) +$$

$$+ \int_0^{x_1} (x_1 - \alpha) D_1^2 u(\alpha, 0) d\alpha + x_2 \int_0^{x_1} (x_1 - \alpha) D_1^2 D_2 u(\alpha, 0) d\alpha +$$

$$+ \int_0^{x_2} (x_2 - \beta) D_2^2 u(0, \beta) d\beta + x_1 \int_0^{x_2} (x_2 - \beta) D_1 D_2^2 u(0, \beta) d\beta +$$

$$+ \int_0^{x_1} \int_0^{x_2} (x_1 - \alpha)(x_2 - \beta) D_1^2 D_2^2 u(\alpha, \beta) d\alpha d\beta.$$